\documentclass[11pt]{article}
\usepackage{latexsym}
\usepackage{amsfonts}
\usepackage{mathrsfs,amsthm,amsmath}
\usepackage{color}
\usepackage{todonotes}
\usepackage{comment}
\newtheorem{theorem}{Theorem}[section]
\newtheorem{proposition}[theorem]{Proposition}

\newtheorem{lemma}{Lemma}[section]

\newtheorem{remark}{Remark}[section]
\newtheorem{example}{Example}[section]
\newtheorem{condition}{Condition}[section]
\begin{document}
\title{Asymptotic results for the absorption time of telegraph processes with a non-standard barrier at the origin}
\author{Antonella Iuliano\thanks{Address: Dipartimento di Matematica, Informatica ed Economia, Universit\`a degli 
Studi della Basilicata, Via dell'Ateneo Lucano 10 (Campus di Macchia Romana), I-85100 Potenza, Italy. e-mail: 
\texttt{antonella.iuliano@unibas.it}}
\and Claudio Macci\thanks{Address: Dipartimento di Matematica, Universit\`a di Roma Tor Vergata, Via della Ricerca 
Scientifica, I-00133 Rome, Italy. e-mail: \texttt{macci@mat.uniroma2.it}}}
\maketitle
\begin{abstract}
A telegraph process with an elastic barrier at the origin was studied in \cite{DicrescenzoMartinucciZacks}; in particular
the number of visits of the origin before the absorption is a geometric distributed random variable $M$.
Some asymptotic results (large and moderate deviations) for that model were obtained in \cite{MacciMartinucciPirozzi}.
In this paper we study large and moderate deviations for a generalized model where $M$ is a light-tailed distributed
random variable.
\end{abstract}

\section{Introduction}
The (integrated) telegraph process describes an alternating random motion of a particle on the real line with finite velocity;
see e.g. the seminal papers \cite{Goldstein} and \cite{Kac}, and the quite recent books \cite{KolesnikRatanov} and \cite{Zacks}.
This model has been widely studied, and it is possible to find several non-standard versions in the literature (generalizations, 
modifications, etc.). Among the more recent references with some generalizations, we recall \cite{CrimaldiDicrescenzoIulianoMartinucci}
for a model driven by certain random trials, \cite{DicrescenzoZacks2015} for a telegraph process perturbed by a Brownian motion, and 
\cite{GarraOrsingher2014} for certain multivariate extensions. Finally, since in this paper we prove results on large deviations, we 
also recall \cite{Macci} (see also some references cited therein). This kind of processes deserves interest for many possible 
applications in different fields.

A large class of (possibly non-standard) telegraph processes concern random motions on the real line subject to barriers; see 
e.g. \cite{MasoliverPorraWeiss}, \cite{FoongKanno}, \cite{Orsingher}  (for results and other references the interested reader 
can see Chapter 3 in \cite{KolesnikRatanov}). There is a wide literature on stochastic processes subject to the presence of barriers 
of different kinds. In particular in some references the barriers exhibit a hard reflection, with random switching to full absorption, 
and they are called \emph{elastic barriers} (see e.g. the old paper \cite{Feller} and the book \cite{Bharucha-Reid}); in this paper we 
use this term even if nowadays it is used for different purposes. In general the number of visits of an elastic barrier is a geometric 
distributed random variable $M$ (say), independent of all the rest. Among the references on stochastic processes subject to elastic 
barriers, here we recall \cite{DicrescenzoMartinucciZacks} and \cite{DicrescenzoMartinucciParaggioZacks} (which deal with some versions 
of telegraph processes), and \cite{GiornoNobilePirozziRicciardi}, \cite{Domine1995} and \cite{Domine1996} (which deal with diffusion 
processes).

In this paper we consider a generalization of the random motion in \cite{DicrescenzoMartinucciZacks}, i.e. a telegraph process 
on $[0,\infty)$, which starts at $x>0$, and with an elastic barrier at the origin. The dynamic of this model depends on two parameters
$\lambda,\mu>0$ such that $\lambda>\mu$. More precisely we consider a non-standard barrier at the origin such that the number of visits
before the absorption is a light-tailed distributed random variable $M$ (thus, in particular, $M$ could be geometric distributed as happens 
for the case of elastic barriers); we are not aware of any other work with
models having this kind of barriers. Our aim is to generalize
the large (and moderate) deviation results in \cite{MacciMartinucciPirozzi} for the model in \cite{DicrescenzoMartinucciZacks}, and more
precisely for the absorption time at the origin $A_x(\lambda,\mu)$. The theory of large deviations provides a collection of techniques 
which allow to give an asymptotic evaluation of the probabilities of rare events on an exponential scale (see e.g. \cite{DemboZeitouni} 
as a reference on this topic). The asymptotic results concern two scalings:
\begin{itemize}
	\item Scaling 1: $x\to\infty$ for $\frac{A_x(\lambda,\mu)}{x}$;
	\item Scaling 2: $\mu\to\infty$ for $A_x(\beta\mu,\mu)$ (for some $\beta>1$).
\end{itemize}

We conclude with the outline of the paper. We start with Section \ref{sec:preliminaries} in which we present some preliminaries on the 
model, some examples, some preliminaries on large deviations and a brief description of the results. In Sections \ref{sec:scaling-1} and 
\ref{sec:scaling-2} we present the results for scalings 1 and 2, respectively. Finally, in Section \ref{sec:numerical-results}, we 
present some numerical estimates based on an asymptotic Normality result under the scaling 2 when $M$ is a shifted Poisson distributed 
random variable (see Example \ref{ex:shifted-Poisson-distribution}).

\section{Preliminaries}\label{sec:preliminaries}
In this section we present the model (with some useful related formulas) and some examples; moreover we present some preliminaries on 
large deviations, and a brief description of the results.

\subsection{The model, some formulas and examples}
Let $\lambda,\mu>0$ be such that $\lambda>\mu$. We consider a random motion of a particle that starts at $x>0$, moves on $[0,\infty)$, 
and we are interested in the absorption time $A_x=A_x(\lambda,\mu)$ at the origin. We refer to eq. (2) in \cite{DicrescenzoMartinucciZacks}
(even if here the distribution of the random variable $M$ could be more general than the one in \cite{DicrescenzoMartinucciZacks}), and 
the absorption time can expressed as follows
$$A_x=C_x+1_{\{M>1\}}\sum_{i=1}^{M-1}C_{0,i},$$
where $C_x,M,\{C_{0,i}:i\geq 1\}$ are independent random variables. More precisely $C_x$ is the random time until the first arrival at 
the origin, and $\{C_{0,i}:i\geq 1\}$ are i.i.d. random variables such that, for every integer $i\geq 1$, $C_{0,i}$ is the (possible) 
$i$-th interarrival time between two consecutive visits of the origin after the time $C_x$; moreover the particle is absorbed at the 
origin after $M$ visits of the origin. In view of what follows it is useful to recall that the moment generating function of $C_x$ is
\begin{equation}\label{eq:mgf-Cx-and-C0}
G_{C_x(\lambda,\mu)}(s)=G_{C_0(\lambda,\mu)}(s)e^{x\Lambda(s;\lambda,\mu)}\ (\mbox{for all}\ s\in\mathbb{R}),
\end{equation}
where
$$G_{C_0(\lambda,\mu)}(s)=\left\{\begin{array}{ll}
\frac{\lambda+\mu-2s-\sqrt{(\lambda+\mu-2s)^2-4\lambda\mu}}{2\mu}&\ \mbox{if}\ s\leq\frac{(\sqrt{\lambda}-\sqrt{\mu})^2}{2}\\
\infty&\ \mbox{if}\ s>\frac{(\sqrt{\lambda}-\sqrt{\mu})^2}{2}
\end{array}\right.$$
is the (common) moment generating function of the random variables $\{C_{0,i}:i\geq 1\}$, and
\begin{equation}\label{eq:def-Lambda}
\Lambda(s;\lambda,\mu):=\frac{\lambda-\mu-\sqrt{(\lambda+\mu-2s)^2-4\lambda\mu}}{2}
\end{equation}
(see eqs. (25)-(41) in \cite{DicrescenzoMartinucciZacks}; see also eq. (3) in \cite{MacciMartinucciPirozzi}).
In what follows we use the symbols $\Lambda^\prime(s;\lambda,\mu)$ and $\Lambda^{\prime\prime}(s;\lambda,\mu)$
for the first and the second derivative of $\Lambda(s;\lambda,\mu)$ with respect to $s$.

In this paper $M$ is a quite general light-tailed random variable according to the following Condition \ref{cond:*}.
Such a condition allows to generalize the case studied in the literature (see eq. (1) in \cite{DicrescenzoMartinucciZacks}) 
where $M$ is a geometric distributed random variable; for more details on this case see also Example 
\ref{ex:geometric-distribution} presented below.

\begin{condition}\label{cond:*}
Let $M$ be a positive and integer valued random variable and assume that there exists $s>0$ such that
$$s\in\mathcal{D}(G_M):=\{r\in\mathbb{R}:G_M(r)<\infty\}.$$
\end{condition}

\begin{remark}\label{rem:*}
Note that in general $G_M$ is an increasing function such that $G_M(0)=1$, and therefore we have 
$(-\infty,0]\subset\mathcal{D}(G_M)$. Then, if Condition \ref{cond:*} holds, we have
$$s_M:=\sup\mathcal{D}(G_M)>0\ (\mbox{possibly with}\ s_M=\infty).$$
Moreover we can have $\mathcal{D}(G_M)=(-\infty,s_M)$, possibly with $s_M=\infty$, or $\mathcal{D}(G_M)=(-\infty,s_M]$
with $s_M<\infty$. For the first case see Example \ref{ex:shifted-Poisson-distribution} (where $s_M=\infty$) and 
Example \ref{ex:geometric-distribution} where $s_M=-\log(1-\alpha)$ for $\alpha\in(0,1]$ (and therefore $s_M=\infty$ 
if and only if $\alpha=1$); for the second case see Example \ref{ex:shifted-Poisson-inverse-Gaussian-distribution}
where $s_M=\log\left(1+\frac{\xi^2}{2\theta}\right)$ for some $\theta,\xi>0$.
\end{remark}

Now we present a formula for the moment generating function $G_{A_x(\lambda,\mu)}$ of $A_x=A_x(\lambda,\mu)$.

\begin{proposition}\label{prop:mgf}
Assume that Condition \ref{cond:*} holds. Then we have
$$G_{A_x(\lambda,\mu)}(s)=\left\{\begin{array}{ll}
G_M(\log G_{C_0(\lambda,\mu)}(s))e^{x\Lambda(s;\lambda,\mu)}&\ \mbox{if}\ s\leq\frac{(\sqrt{\lambda}-\sqrt{\mu})^2}{2}\\
\infty&\ \mbox{if}\ s>\frac{(\sqrt{\lambda}-\sqrt{\mu})^2}{2}.
\end{array}\right.$$
\end{proposition}
\begin{proof}
Firstly we have
$$G_{A_x(\lambda,\mu)}(s)=\mathbb{E}\left[e^{s(C_x+1_{\{M>1\}}\sum_{i=1}^{M-1}C_{0,i})}\right]=
G_{C_x(\lambda,\mu)}(s)\mathbb{E}\left[e^{s1_{\{M>1\}}\sum_{i=1}^{M-1}C_{0,i}}\right]$$
and, by \eqref{eq:mgf-Cx-and-C0}, we get
$$G_{A_x(\lambda,\mu)}(s)=G_{C_0(\lambda,\mu)}(s)e^{x\Lambda(s;\lambda,\mu)}\mathbb{E}\left[e^{s1_{\{M>1\}}\sum_{i=1}^{M-1}C_{0,i}}\right].$$
Moreover we have
\begin{multline*}
\mathbb{E}\left[e^{s1_{\{M>1\}}\sum_{i=1}^{M-1}C_{0,i}}\right]
=\mathbb{E}\left[\mathbb{E}\left[e^{s1_{\{M>1\}}\sum_{i=1}^{M-1}C_{0,i}}|M\right]\right]\\
=P(M=1)+\sum_{m=2}^\infty(G_{C_0(\lambda,\mu)}(s))^{m-1}P(M=m)=\sum_{m=1}^\infty(G_{C_0(\lambda,\mu)}(s))^{m-1}P(M=m)\\
=\frac{\sum_{m=1}^\infty(G_{C_0(\lambda,\mu)}(s))^mP(M=m)}{G_{C_0(\lambda,\mu)}(s)}
=\frac{G_M(\log G_{C_0(\lambda,\mu)}(s))}{G_{C_0(\lambda,\mu)}(s)}.
\end{multline*}
Then we conclude by combining the above equalities.
\end{proof}

It is worth noting that we can have $G_{A_x(\lambda,\mu)}(s)=\infty$ for some $s\in\left(0,\frac{(\sqrt{\lambda}-\sqrt{\mu})^2}{2}\right]$.
This issue is discussed in the next remark.

\begin{remark}\label{rem:the-domain-of-mgf-of-absorption-time}
Assume that Condition \ref{cond:*} holds. Then, if we consider the set
$$\mathcal{D}(G_{A_x(\lambda,\mu)}):=\{r\in\mathbb{R}:G_{A_x(\lambda,\mu)}(r)<\infty\},$$
we have
$$(-\infty,0]\subset\mathcal{D}(G_{A_x(\lambda,\mu)})\subset\left(-\infty,\frac{(\sqrt{\lambda}-\sqrt{\mu})^2}{2}\right],$$
and the first inclusion is strict. Moreover $\log G_{C_0(\lambda,\mu)}$ is an increasing function; then we can consider its inverse 
$[\log G_{C_0(\lambda,\mu)}]^{-1}(s)$. In particular we have
$$G_{C_0(\lambda,\mu)}\left(\frac{(\sqrt{\lambda}-\sqrt{\mu})^2}{2}\right)=\sqrt{\frac{\lambda}{\mu}},
\ \mbox{which yields}\ [\log G_{C_0(\lambda,\mu)}]^{-1}\left(\log\sqrt{\frac{\lambda}{\mu}}\right)=\frac{(\sqrt{\lambda}-\sqrt{\mu})^2}{2}.$$
Then, by Proposition \ref{prop:mgf}, we have two cases.
\begin{itemize}
\item Case A. If $s_M\geq\log\sqrt{\frac{\lambda}{\mu}}$, possibly with $s_M=\infty$, we have:
$$\mbox{if $s_M=\log\sqrt{\frac{\lambda}{\mu}}$ and $\mathcal{D}(G_M)=(-\infty,s_M)$, then 
$\mathcal{D}(G_{A_x(\lambda,\mu)})=\left(-\infty,\frac{(\sqrt{\lambda}-\sqrt{\mu})^2}{2}\right)$};$$
otherwise
$$\mathcal{D}(G_{A_x(\lambda,\mu)})=\left(-\infty,\frac{(\sqrt{\lambda}-\sqrt{\mu})^2}{2}\right].$$
\item Case B. If $s_M<\log\sqrt{\frac{\lambda}{\mu}}$, then we set
\begin{equation}\label{eq:def-s-hat}
\hat{s}(\lambda,\mu,s_M):=[\log G_{C_0(\lambda,\mu)}]^{-1}(s_M)
\end{equation}
and we have
$$\mathcal{D}(G_{A_x(\lambda,\mu)})=\left\{\begin{array}{ll}
(-\infty,\hat{s}(\lambda,\mu,s_M))&\ \mbox{if}\ \mathcal{D}(G_M)=(-\infty,s_M)\\
(-\infty,\hat{s}(\lambda,\mu,s_M)]&\ \mbox{if}\ \mathcal{D}(G_M)=(-\infty,s_M].
\end{array}\right.$$
\end{itemize}
\end{remark}

Now we present three examples.

\begin{example}\label{ex:shifted-Poisson-distribution}
Assume that $M$ is a shifted Poisson distributed random variable; thus, for some $\theta>0$,
$$P(M=m)=\frac{\theta^{m-1}}{(m-1)!}e^{-\theta}\ (\mbox{for all}\ m\geq 1).$$
Then we can easily check that
$$G_M(s)=e^{s+\theta(e^s-1)}\ \mbox{for all}\ s\in\mathbb{R};$$
thus $s_M=\infty$ and $\mathcal{D}(G_M)=\mathbb{R}$ (so this example obviously concerns the Case A in Remark 
\ref{rem:the-domain-of-mgf-of-absorption-time}). Then, by Proposition \ref{prop:mgf}, we can easily check that
$$G_{A_x(\lambda,\mu)}(s)=\left\{\begin{array}{ll}
G_{C_0(\lambda,\mu)}(s)e^{\theta(G_{C_0(\lambda,\mu)}(s)-1)}e^{x\Lambda(s;\lambda,\mu)}&\ \mbox{if}\ s\leq\frac{(\sqrt{\lambda}-\sqrt{\mu})^2}{2}\\
\infty&\ \mbox{if}\ s>\frac{(\sqrt{\lambda}-\sqrt{\mu})^2}{2}.
\end{array}\right.$$
\end{example}

\begin{example}\label{ex:shifted-Poisson-inverse-Gaussian-distribution}
Assume that $M$ is a shifted Poisson inverse Gaussian distributed random variable. Then, for some $\theta,\xi>0$,
the moment generating function of $M$ is
$$G_M(s)=\left\{\begin{array}{ll}
e^{s+\xi-\sqrt{\xi^2-2\theta(e^s-1)}}&\mbox{if}\ \theta(e^s-1)\leq\frac{\xi^2}{2},\ \mbox{i.e.}\ s\leq\log\left(1+\frac{\xi^2}{2\theta}\right)\\
\infty& \mbox{if}\ \theta(e^s-1)>\frac{\xi^2}{2},\ \mbox{i.e.}\ s>\log\left(1+\frac{\xi^2}{2\theta}\right);
\end{array}\right.$$
thus $\mathcal{D}(G_M)=(-\infty,s_M]$ with $s_M=\log\left(1+\frac{\xi^2}{2\theta}\right)$. As far as Case A in Remark 
\ref{rem:the-domain-of-mgf-of-absorption-time} is concerned, we have $s_M\geq\log\sqrt{\frac{\lambda}{\mu}}$ if and 
only if $1+\frac{\xi^2}{2\theta}\geq\sqrt{\frac{\lambda}{\mu}}$ or, equivalently, if and only if 
$\frac{\xi^2}{2\theta}\geq\sqrt{\frac{\lambda}{\mu}}-1$. Then, by Proposition \ref{prop:mgf}, we can easily check 
that we have the following two cases (which correspond to Cases A and B in Remark \ref{rem:the-domain-of-mgf-of-absorption-time},
respectively).
\begin{itemize}
\item If $\frac{\xi^2}{2\theta}\geq\sqrt{\frac{\lambda}{\mu}}-1$, then
$$G_{A_x(\lambda,\mu)}(s)=\left\{\begin{array}{ll}
G_{C_0(\lambda,\mu)}(s)e^{\xi-\sqrt{\xi^2-2\theta(G_{C_0(\lambda,\mu)}(s)-1)}+x\Lambda(s;\lambda,\mu)}
&\ \mbox{if}\ s\leq\frac{(\sqrt{\lambda}-\sqrt{\mu})^2}{2}\\
\infty&\ \mbox{otherwise}.
\end{array}\right.$$
\item If $\frac{\xi^2}{2\theta}<\sqrt{\frac{\lambda}{\mu}}-1$, then
$$G_{A_x(\lambda,\mu)}(s)=\left\{\begin{array}{ll}
G_{C_0(\lambda,\mu)}(s)e^{\xi-\sqrt{\xi^2-2\theta(G_{C_0(\lambda,\mu)}(s)-1)}+x\Lambda(s;\lambda,\mu)}
&\ \mbox{if}\ s\leq\hat{s}(\lambda,\mu,s_M)\\
\infty&\ \mbox{otherwise}.
\end{array}\right.$$
Note that (see Case B in Remark \ref{rem:the-domain-of-mgf-of-absorption-time}), since
$\mathcal{D}(G_M)=(-\infty,s_M]$ with $s_M=\log\left(1+\frac{\xi^2}{2\theta}\right)$, we have
$\mathcal{D}(G_{A_x(\lambda,\mu)})=(-\infty,\hat{s}(\lambda,\mu,s_M)]$. In particular we can
check that
$$G_{A_x(\lambda,\mu)}(\hat{s}(\lambda,\mu,s_M))=\left(1+\frac{\xi^2}{2\theta}\right)
e^{\xi+x\Lambda(\hat{s}(\lambda,\mu,s_M);\lambda,\mu)}<\infty$$
since $G_{C_0(\lambda,\mu)}(\hat{s}(\lambda,\mu,s_M))=1+\frac{\xi^2}{2\theta}$.
\end{itemize}
\end{example}

\begin{example}\label{ex:geometric-distribution}
Assume that $M$ is a geometric distributed random variable; thus, for some $\alpha\in(0,1]$,
$$P(M=m)=(1-\alpha)^{m-1}\alpha\ (\mbox{for all}\ m\geq 1)$$
as in eq. (1) in \cite{DicrescenzoMartinucciZacks}. Then we can easily compute the moment generating function 
of $M$, and we have
$$G_M(s)=\left\{\begin{array}{ll}
\frac{\alpha e^s}{1-(1-\alpha)e^s}&\mbox{if}\ (1-\alpha)e^s<1,\ \mbox{i.e.}\ s<-\log(1-\alpha)\\
\infty& \mbox{if}\ (1-\alpha)e^s\geq 1,\ \mbox{i.e.}\ s\geq-\log(1-\alpha);
\end{array}\right.$$
thus $\mathcal{D}(G_M)=(-\infty,s_M)$ with $s_M=-\log(1-\alpha)$, and we have $s_M=\infty$ if and only if $\alpha=1$ 
(because we consider the rule $\log 0=-\infty$). As far as Case A in Remark \ref{rem:the-domain-of-mgf-of-absorption-time}
is concerned, we have $s_M\geq\log\sqrt{\frac{\lambda}{\mu}}$ if and only if $\frac{1}{1-\alpha}\geq\sqrt{\frac{\lambda}{\mu}}$ or,
equivalently, if and only if $\alpha\geq 1-\sqrt{\frac{\mu}{\lambda}}$. Then, by Proposition \ref{prop:mgf}, we can easily check 
that we have the following two cases (which correspond to Cases A and B in Remark \ref{rem:the-domain-of-mgf-of-absorption-time},
respectively).
\begin{itemize}
\item If $\alpha\geq 1-\sqrt{\frac{\mu}{\lambda}}$, then
$$G_{A_x(\lambda,\mu)}(s)=\left\{\begin{array}{ll}
\frac{\alpha G_{C_0(\lambda,\mu)}(s)e^{x\Lambda(s;\lambda,\mu)}}{1-(1-\alpha)G_{C_0(\lambda,\mu)}(s)}
&\ \mbox{if}\ s\leq\frac{(\sqrt{\lambda}-\sqrt{\mu})^2}{2}\\
\infty&\ \mbox{otherwise}.
\end{array}\right.$$
In particular, if $\alpha=1-\sqrt{\frac{\mu}{\lambda}}$, we can write down
$$G_{A_x(\lambda,\mu)}(s)=\left\{\begin{array}{ll}
\frac{\alpha G_{C_0(\lambda,\mu)}(s)e^{x\Lambda(s;\lambda,\mu)}}{1-(1-\alpha)G_{C_0(\lambda,\mu)}(s)}
&\ \mbox{if}\ s<\frac{(\sqrt{\lambda}-\sqrt{\mu})^2}{2}\\
\infty&\ \mbox{otherwise}
\end{array}\right.$$
because $G_{A_x(\lambda,\mu)}\left(\frac{(\sqrt{\lambda}-\sqrt{\mu})^2}{2}\right)=\infty$, and therefore
$\mathcal{D}(G_{A_x(\lambda,\mu)})=\left(-\infty,\frac{(\sqrt{\lambda}-\sqrt{\mu})^2}{2}\right)$. 
\item If $\alpha<1-\sqrt{\frac{\mu}{\lambda}}$, then
\begin{multline*}
G_{A_x(\lambda,\mu)}(s)=\left\{\begin{array}{ll}
\frac{\alpha G_{C_0(\lambda,\mu)}(s)e^{x\Lambda(s;\lambda,\mu)}}{1-(1-\alpha)G_{C_0(\lambda,\mu)}(s)}
&\ \mbox{if}\ s\leq\hat{s}(\lambda,\mu,s_M)\\
\infty&\ \mbox{otherwise}
\end{array}\right.\nonumber\\
=\left\{\begin{array}{ll}
\frac{\alpha G_{C_0(\lambda,\mu)}(s)e^{x\Lambda(s;\lambda,\mu)}}{1-(1-\alpha)G_{C_0(\lambda,\mu)}(s)}
&\ \mbox{if}\ s<\hat{s}(\lambda,\mu,s_M),\ \mbox{i.e.}\ (1-\alpha)G_{C_0(\lambda,\mu)}(s)<1\\
\infty&\ \mbox{otherwise},
\end{array}\right.
\end{multline*}
where the second equality holds noting that 
$\mathcal{D}(G_{A_x(\lambda,\mu)})=(-\infty,\hat{s}(\lambda,\mu,s_M))$ because 
$\mathcal{D}(G_M)=(-\infty,s_M)$ (see Case B in Remark \ref{rem:the-domain-of-mgf-of-absorption-time})
and, moreover, since $s_M=-\log(1-\alpha)$, we have $s<\hat{s}(\lambda,\mu,s_M)$ if and only if 
$G_{C_0(\lambda,\mu)}(s)<e^{s_M}=\frac{1}{1-\alpha}$. Finally we remark that the equality 
$G_{A_x(\lambda,\mu)}(\hat{s}(\lambda,\mu,s_M))=\infty$ can be checked noting that
$\alpha G_{C_0(\lambda,\mu)}(\hat{s}(\lambda,\mu,s_M))=\frac{\alpha}{1-\alpha}$ and
$1-(1-\alpha)G_{C_0(\lambda,\mu)}(\hat{s}(\lambda,\mu,s_M))=0$.
\end{itemize}
Note that, in both cases $\alpha\geq 1-\sqrt{\frac{\mu}{\lambda}}$ and $\alpha<1-\sqrt{\frac{\mu}{\lambda}}$, for some
values of $s$ we have
$$G_{A_x(\lambda,\mu)}(s)=\frac{\alpha G_{C_0(\lambda,\mu)}(s)e^{x\Lambda(s;\lambda,\mu)}}{1-(1-\alpha)G_{C_0(\lambda,\mu)}(s)}
=\frac{\alpha G_{C_x(\lambda,\mu)}(s)}{1+(\alpha-1)G_{C_0(\lambda,\mu)}(s)}$$
(for the second equality see eq. \eqref{eq:mgf-Cx-and-C0}). In this way we recover the first displayed formula in the proof
of Proposition 9 in \cite{DicrescenzoMartinucciZacks}. Actually one should consider the correct version of Proposition 9 in 
\cite{DicrescenzoMartinucciZacks} discussed in Remark 2.1 in \cite{MacciMartinucciPirozzi}; in particular, for the case 
$\alpha<1-\sqrt{\frac{\mu}{\lambda}}$, $\hat{s}(\lambda,\mu,\alpha)=\frac{\alpha((1-\alpha)\lambda-\mu)}{2(1-\alpha)}$ 
in Remark 2.1 in \cite{MacciMartinucciPirozzi} corresponds to $\hat{s}(\lambda,\mu,s_M)$ in Remark 
\ref{rem:the-domain-of-mgf-of-absorption-time} (Case B) in this paper.
\end{example}

\subsection{Preliminaries on large deviations, and a brief description of the results}
We start with some basic definitions (see e.g. \cite{DemboZeitouni}, pages 4-5). Let $\mathcal{Z}$ be a topological space 
equipped with its completed Borel $\sigma$-field. A family of $\mathcal{Z}$-valued random variables $\{Z_r:r>0\}$ (defined 
on the same probability space $(\Omega,\mathcal{F},P)$) satisfies the large deviation principle (LDP for short) with speed 
function $v_r$ and rate function $I$ if: $\lim_{r\to\infty}v_r=\infty$; the function $I:\mathcal{Z}\to[0,\infty]$ is lower 
semi-continuous;
\begin{equation}\label{eq:UB-LDP-definition}
\limsup_{n\to\infty}\frac{1}{v_r}\log P(Z_r\in F)\leq-\inf_{z\in F}I(z)\ \mbox{for all closed sets}\ F
\end{equation}
and
\begin{equation}\label{eq:LB-LDP-definition}
\liminf_{r\to\infty}\frac{1}{v_r}\log P(Z_r\in G)\geq-\inf_{z\in G}I(z)\ \mbox{for all open sets}\ G.
\end{equation}
A rate function $I$ is said to be \emph{good} if its level sets
$\{\{z\in\mathcal{Z}:I(z)\leq\eta\}:\eta\geq 0\}$ are compact.

Throughout this paper we prove LDPs with $\mathcal{Z}=\mathbb{R}$. In view of what follows we recall a well-known result
(specified for real-valued random variables) which provides \eqref{eq:UB-LDP-definition} and a weak form of 
\eqref{eq:LB-LDP-definition} with $I=\Lambda^*$.

\begin{theorem}[G\"{a}rtner Ellis Theorem (on $\mathbb{R}$)]\label{th:GE}
Let $\{Z_r:r>0\}$ be a family of real valued random variables (defined on the same probability space $(\Omega,\mathcal{F},P)$).
Assume that the function $\Lambda:\mathbb{R}\to(-\infty,\infty]$ defined by
$$\Lambda(s):=\lim_{r\to\infty}\frac{1}{v_r}\log\mathbb{E}\left[e^{v_rs Z_r}\right]\ (\mbox{for all}\ s\in\mathbb{R})$$
exists, and it is finite in a neighborhood of the origin $s=0$. Moreover let $\Lambda^*:\mathbb{R}\to[0,\infty]$ defined by
\begin{equation}\label{eq:LegendreTransform}
\Lambda^*(z):=\sup_{s\in\mathbb{R}}\{sz-\Lambda(s)\}
\end{equation}
(it is the Legendre transform of $\Lambda$). Then: \eqref{eq:UB-LDP-definition} holds with $I=\Lambda^*$;
$$\liminf_{r\to\infty}\frac{1}{v_r}\log P(Z_r\in G)\geq-\inf_{z\in G\cap\mathcal{E}}\Lambda^*(z)\ \mbox{for all open sets}\ G$$
where $\mathcal{E}$ is the set of exposed points of $I$ (namely the points in which $I$ is finite and strictly convex); if 
$\Lambda$ is essentially smooth and lower semi-continuous, then the LDP holds with good rate function $I=\Lambda^*$.
\end{theorem}

We also recall that $\Lambda$ in the above statement is essentially smooth (see e.g. Definition 2.3.5 in \cite{DemboZeitouni}) if:
\begin{itemize}
\item the interior of the set $\mathcal{D}_\Lambda:=\{s\in\mathbb{R}:\Lambda(s)<\infty\}$ is non-empty;
\item the function $\Lambda$ is differentiable throughout the interior of $\mathcal{D}_\Lambda$;
\item the function $\Lambda$  is a steep (namely $|\Lambda^\prime(s)|$ tends to infinity when $s$ in the interior of 
$\mathcal{D}_\Lambda$ approaches any finite point of its boundary).
\end{itemize}

For instance the function $\Lambda(\cdot;\lambda,\mu)$ in \eqref{eq:def-Lambda} is essentially smooth because 
$\Lambda^\prime(s;\lambda,\mu)\uparrow\infty$ as $s\uparrow\frac{(\sqrt{\lambda}-\sqrt{\mu})^2}{2}$ (this can be checked with some 
computations and we omit the details).

Now, in view of what follows, we present some formulas for Legendre transforms (see \eqref{eq:LegendreTransform}). This is the analogue
of Lemma 2.1 in \cite{MacciMartinucciPirozzi} (in some parts we have exactly the same formulas, in other cases some notation are suitably
changed) and we omit the proof. Note that the two cases presented in the following lemma correspond to Cases A and B in Remark \ref{rem:the-domain-of-mgf-of-absorption-time}, respectively.

\begin{lemma}\label{lem:Legendre}
Let $\Lambda(\cdot;\lambda,\mu)$ be the function in \eqref{eq:def-Lambda}.\\
(i) Let $H_A(z;\lambda,\mu)$ be defined by 
$$H_A(z;\lambda,\mu):=\sup_{s\leq\frac{(\sqrt{\lambda}-\sqrt{\mu})^2}{2}}\{sz-\Lambda(s;\lambda,\mu)\},$$
or equivalently
$$H_A(z;\lambda,\mu):=\sup_{s<\frac{(\sqrt{\lambda}-\sqrt{\mu})^2}{2}}\{sz-\Lambda(s;\lambda,\mu)\}.$$
Then we have
$$H_A(z;\lambda,\mu)=\left\{\begin{array}{ll}
\frac{1}{2}\left(\sqrt{(z-1)\lambda}-\sqrt{(z+1)\mu}\right)^2&\ \mbox{if}\ z\geq 1\\
\infty&\ \mbox{otherwise}.
\end{array}\right.$$
(ii) For $s_M<\log\sqrt{\frac{\lambda}{\mu}}$, let $\hat{s}(\lambda,\mu,s_M)$ be defined by \eqref{eq:def-s-hat}
(see Case B in Remark \ref{rem:the-domain-of-mgf-of-absorption-time}), and set
$$\tilde{z}(\lambda,\mu,s_M):=\Lambda^\prime(\hat{s}(\lambda,\mu,s_M);\lambda,\mu).$$
Moreover let $H_B(z;\lambda,\mu,s_M)$ be defined by 
$$H_B(z;\lambda,\mu,s_M):=\sup_{s\leq\hat{s}(\lambda,\mu,s_M)}\{sz-\Lambda(s;\lambda,\mu)\},$$
or equivalently
$$H_B(z;\lambda,\mu,s_M):=\sup_{s<\hat{s}(\lambda,\mu,s_M)}\{sz-\Lambda(s;\lambda,\mu)\}.$$
Then we have
\begin{multline*}
H_B(z;\lambda,\mu,s_M)=\left\{\begin{array}{ll}
H_A(z;\lambda,\mu)&\ \mbox{if}\ z\leq\tilde{z}(\lambda,\mu,s_M)\\
\hat{s}(\lambda,\mu,s_M)z-\Lambda(\hat{s}(\lambda,\mu,s_M);\lambda,\mu)&\ \mbox{otherwise}
\end{array}\right.\\
=\left\{\begin{array}{ll}
\infty&\ \mbox{if}\ z<1\\
\frac{1}{2}\left(\sqrt{(z-1)\lambda}-\sqrt{(z+1)\mu}\right)^2&\ \mbox{if}\ 1\leq z\leq\tilde{z}(\lambda,\mu,s_M)\\
\hat{s}(\lambda,\mu,s_M)z-\Lambda(\hat{s}(\lambda,\mu,s_M);\lambda,\mu)&\ \mbox{if}\ z>\tilde{z}(\lambda,\mu,s_M).
\end{array}\right.
\end{multline*}
\end{lemma}

\begin{remark}\label{rem:geometric-distribution-part-2}
One can check with some computations that
$$\tilde{z}(\lambda,\mu,s_M):=\frac{\lambda+\mu-2\hat{s}(\lambda,\mu,s_M)}
{\sqrt{(\lambda+\mu-2\hat{s}(\lambda,\mu,s_M))^2-4\lambda\mu}},$$
which is the analogue of $\tilde{z}(\lambda,\mu,\alpha)$ in Lemma 2.1 in \cite{MacciMartinucciPirozzi}.
\end{remark}

Our aim is to present a generalization of some large and moderate deviation results in \cite{MacciMartinucciPirozzi}.
These asymptotic results concern the random variables $\left\{\frac{A_x(\lambda,\mu)}{x}:x>0\right\}$ as the initial 
position $x$ go to infinity (scaling 1), and the random variables $\{A_x(\beta\mu,\mu):\mu>0\}$ as the switching rates 
$\lambda=\beta\mu$ and $\mu$ go to infinity simultaneously (scaling 2).

Lemma \ref{lem:Legendre} will be used in Proposition \ref{prop:LD1} and Remark \ref{rem:LD1} for scaling 1 
(Cases A and B respectively) with speed $x$, and in Proposition \ref{prop:LD2} and Remark \ref{rem:LD2}
for scaling 2 (Cases A and B respectively) with speed $\mu$. Some other results concern \emph{moderate 
deviations} for which we do not have to distinguish between Cases A and B: Proposition \ref{prop:MD1}
for scaling 1, and Proposition \ref{prop:MD2} for scaling 2. Actually, for scaling 2, we also present a
non-central moderate deviation result, i.e. Proposition \ref{prop:MD2-nc} together with Lemma 
\ref{lem:weak-convergence}.

The moderate deviation results in Propositions \ref{prop:MD1} and \ref{prop:MD2} fill the gap between a 
convergence to a constant (governed by a suitable LDP) and a weak convergence to a Normal distribution; 
for more details see Remarks 3.1 and 4.2 in \cite{MacciMartinucciPirozzi}, respectively. We can also say
that the non-central moderate deviation result in Proposition \ref{prop:MD2-nc} fill the gap between a
convergence to a constant (governed by a suitable LDP) and the weak convergence of 
$\mu A_{x/\mu}(\beta\mu,\mu)$ to $A_x(\beta,1)$ as $\mu\to\infty$, which is a consequence of Lemma 
\ref{lem:weak-convergence}.

\section{Large and moderate deviation results under the scaling 1}\label{sec:scaling-1}
We start with the analogue of Proposition 3.1 and Remark 3.2 in \cite{MacciMartinucciPirozzi}; in the
first case we have a full LDP, in the second case we have a weak formulation of the lower bound for 
open sets in term of the exposed points (as illustrated in Theorem \ref{th:GE}).

\begin{proposition}\label{prop:LD1}
Assume that $s_M\geq\log\sqrt{\frac{\lambda}{\mu}}$. Then the family $\left\{\frac{A_x(\lambda,\mu)}{x}:x>0\right\}$ 
satisfies the LDP with speed $x$, and good rate function $I_1$ defined by $I_1(z):=H_A(z;\lambda,\mu)$, where
$H_A(z;\lambda,\mu)$ is the function in Lemma \ref{lem:Legendre}(i).
\end{proposition}
\begin{proof}
We consider Proposition \ref{prop:mgf}, Remark \ref{rem:the-domain-of-mgf-of-absorption-time} (Case A) and Lemma
\ref{lem:Legendre}(i). If $s_M=\log\sqrt{\frac{\lambda}{\mu}}$ and $\mathcal{D}(G_M)$ is open, then
$$\lim_{x\to\infty}\frac{1}{x}\log\mathbb{E}\left[e^{sA_x(\lambda,\mu)}\right]=\left\{\begin{array}{ll}
\Lambda(s;\lambda,\mu)&\ \mbox{if}\ s<\frac{(\sqrt{\lambda}-\sqrt{\mu})^2}{2}\\
\infty&\ \mbox{if}\ s\geq\frac{(\sqrt{\lambda}-\sqrt{\mu})^2}{2};
\end{array}\right.$$
otherwise
$$\lim_{x\to\infty}\frac{1}{x}\log\mathbb{E}\left[e^{sA_x(\lambda,\mu)}\right]=\left\{\begin{array}{ll}
\Lambda(s;\lambda,\mu)&\ \mbox{if}\ s\leq\frac{(\sqrt{\lambda}-\sqrt{\mu})^2}{2}\\
\infty&\ \mbox{if}\ s>\frac{(\sqrt{\lambda}-\sqrt{\mu})^2}{2}.
\end{array}\right.$$
Then the desired LDP holds by a straightforward application of Theorem \ref{th:GE}.
\end{proof}

\begin{remark}\label{rem:LD1}
If $s_M<\log\sqrt{\frac{\lambda}{\mu}}$, then we have to consider Remark \ref{rem:the-domain-of-mgf-of-absorption-time}
(Case B) and the function $H_B(z;\lambda,\mu,s_M)$ in Lemma \ref{lem:Legendre}(ii). Then, by Theorem \ref{th:GE}, we have
$$\limsup_{x\to\infty}\frac{1}{x}\log P\left(\frac{A_x(\lambda,\mu)}{x}\in F\right)
\leq-\inf_{z\in F}H_B(z;\lambda,\mu,s_M)\ \mbox{for all closed sets}\ F$$
and
$$\liminf_{x\to\infty}\frac{1}{x}\log P\left(\frac{A_{x}(\lambda,\mu)}{x}\in G\right)
\geq-\inf_{z\in G\cap\mathcal{E} }H_B(z;\lambda,\mu,s_M)\ \mbox{for all open sets}\ G$$
where $\mathcal{E}=(\tilde{z}(\lambda,\mu,s_M),\infty)$ is the set of exposed points of $H_B(\cdot;\lambda,\mu,s_M)$.
\end{remark}

We conclude with moderate deviations, i.e. with the analogue of Proposition 3.2 in \cite{MacciMartinucciPirozzi}.

\begin{proposition}\label{prop:MD1}
For every family of positive numbers $\{\varepsilon_x:x>0\}$ such that
\begin{equation}\label{eq:MD1-conditions}
\varepsilon_x\to 0\ \mbox{and}\ x\varepsilon_x\to\infty,
\end{equation}
the family $\left\{\frac{A_x(\lambda,\mu)-\mathbb{E}[A_x(\lambda,\mu)]}{\sqrt{x/\varepsilon_x}}:x>0\right\}$ satisfies the LDP with speed
$1/\varepsilon_x$, and good rate function $\tilde{I}_1$ defined by
$\tilde{I}_1(z):=\frac{z^2}{2\Lambda^{\prime\prime}(0;\lambda,\mu)}$,
where $\Lambda^{\prime\prime}(0;\lambda,\mu)=\frac{8\lambda\mu}{(\lambda-\mu)^3}$.
\end{proposition}
\begin{proof}
We follow the lines of the proof of Proposition 3.2 in \cite{MacciMartinucciPirozzi} and we have to prove that
$$\lim_{x\to\infty}\frac{1}{1/\varepsilon_x}
\log\mathbb{E}\left[e^{\frac{s}{\varepsilon_x}\frac{A_x(\lambda,\mu)-\mathbb{E}[A_x(\lambda,\mu)]}{\sqrt{x/\varepsilon_x}}}\right]
=\frac{\Lambda^{\prime\prime}(0;\lambda,\mu)}{2}s^2\ (\mbox{for all}\ s\in\mathbb{R}).$$
Firstly we note that
\begin{multline*}
\frac{1}{1/\varepsilon_x}
\log\mathbb{E}\left[e^{\frac{s}{\varepsilon_x}\frac{A_x(\lambda,\mu)-\mathbb{E}[A_x(\lambda,\mu)]}{\sqrt{x/\varepsilon_x}}}\right]
=\varepsilon_x\left(\log G_{A_x(\lambda,\mu)}\left(\frac{s}{\sqrt{x\varepsilon_x}}\right)
-\mathbb{E}[A_x(\lambda,\mu)]\frac{s}{\sqrt{x\varepsilon_x}}\right)\\
=x\varepsilon_x\left(\frac{1}{x}\log G_{A_x(\lambda,\mu)}\left(\frac{s}{\sqrt{x\varepsilon_x}}\right)-
\frac{\mathbb{E}[A_x(\lambda,\mu)]}{x}\frac{s}{\sqrt{x\varepsilon_x}}\right),
\end{multline*}
where $\frac{s}{\sqrt{x\varepsilon_x}}$ is close to zero if $x$ is large enough. Moreover
\begin{equation}\label{eq:expected-value-1}
\mathbb{E}[A_x(\lambda,\mu)]=[\log G_{A_x(\lambda,\mu)}]^\prime(0)
=G_M^\prime(0)G_{C_0(\lambda,\mu)}^\prime(0)+x\Lambda^\prime(0;\lambda,\mu)
\end{equation}
(the second equality can be checked with some computations). Then, by Proposition \ref{prop:mgf}, 
eq. \eqref{eq:expected-value-1} and the Taylor formula of order 2 for the function $\Lambda(s;\lambda,\mu)$, for $x$ large enough we have
\begin{multline*}
\frac{1}{1/\varepsilon_x}
\log\mathbb{E}\left[e^{\frac{s}{\varepsilon_x}\frac{A_x(\lambda,\mu)-\mathbb{E}[A_x(\lambda,\mu)]}{\sqrt{x/\varepsilon_x}}}\right]\\
=x\varepsilon_x\left\{\frac{1}{x}\log G_M\left(\log G_{C_0(\lambda,\mu)}\left(\frac{s}{\sqrt{x\varepsilon_x}}\right)\right)
+\Lambda\left(\frac{s}{\sqrt{x\varepsilon_x}};\lambda,\mu\right)\right.\\
\left.-\left(\frac{G_M^\prime(0)G_{C_0(\lambda,\mu)}^\prime(0)}{x}+\Lambda^\prime(0;\lambda,\mu)\right)\frac{s}{\sqrt{x\varepsilon_x}}\right\}\\
=\varepsilon_x\left(\log G_M\left(\log G_{C_0(\lambda,\mu)}\left(\frac{s}{\sqrt{x\varepsilon_x}}\right)\right)
-G_M^\prime(0)G_{C_0(\lambda,\mu)}^\prime(0)\frac{s}{\sqrt{x\varepsilon_x}}\right)
+x\varepsilon_x\left(\frac{\Lambda^{\prime\prime}(0;\lambda,\mu)}{2}\frac{s^2}{x\varepsilon_x}+o\left(\frac{1}{x\varepsilon_x}\right)\right)
\end{multline*}
and we conclude by taking the limit as $x\to\infty$.
\end{proof}

\section{Large and moderate deviation results under the scaling 2}\label{sec:scaling-2}
We start with the analogue of Proposition 4.1 and Remark 4.1 in \cite{MacciMartinucciPirozzi}; in the
first case we have a full LDP, in the second case we have a weak formulation of the lower bound for 
open sets in term of the exposed points (as illustrated in Theorem \ref{th:GE}).

\begin{proposition}\label{prop:LD2}
Assume that $s_M\geq\log\sqrt{\beta}$, for $\beta>1$. Then the family $\{A_x(\beta\mu,\mu):\mu>0\}$ 
satisfies the LDP with speed $\mu$, and good rate function $I_2$ defined by $I_2(z):=xH_A(z/x;\beta,1)$, where
$H_A(z;\lambda,\mu)$ is the function in Lemma \ref{lem:Legendre}(i).
\end{proposition}
\begin{proof}
We consider Proposition \ref{prop:mgf} (note that $G_{C_0(\beta\mu,\mu)}(\mu s)=G_{C_0(\beta,1)}(s)$
and $\Lambda(\mu s;\beta\mu,\mu)=\mu\Lambda(s;\beta,1)$), Remark 
\ref{rem:the-domain-of-mgf-of-absorption-time} (Case A) and Lemma \ref{lem:Legendre}(i). If
$s_M=\log\sqrt{\beta}$ and $\mathcal{D}(G_M)$ is open, then
$$\lim_{\mu\to\infty}\frac{1}{\mu}\log\mathbb{E}\left[e^{\mu sA_x(\beta\mu,\mu)}\right]=\left\{\begin{array}{ll}
x\Lambda(s;\beta,1)&\ \mbox{if}\ s<\frac{(\sqrt{\beta}-1)^2}{2}\\
\infty&\ \mbox{if}\ s\geq\frac{(\sqrt{\beta}-1)^2}{2};
\end{array}\right.$$
otherwise
$$\lim_{\mu\to\infty}\frac{1}{\mu}\log\mathbb{E}\left[e^{\mu sA_x(\beta\mu,\mu)}\right]=\left\{\begin{array}{ll}
x\Lambda(s;\beta,1)&\ \mbox{if}\ s\leq\frac{(\sqrt{\beta}-1)^2}{2}\\
\infty&\ \mbox{if}\ s>\frac{(\sqrt{\beta}-1)^2}{2}.
\end{array}\right.$$
Then the desired LDP holds by a straightforward application of Theorem \ref{th:GE}.
\end{proof}

\begin{remark}\label{rem:LD2}
If $s_M<\log\sqrt{\beta}$ for $\beta>1$, then we have to consider Remark \ref{rem:the-domain-of-mgf-of-absorption-time}
(Case B) and the function $H_B(z;\beta,1,s_M)$ in Lemma \ref{lem:Legendre}(ii). Then, by Theorem \ref{th:GE}, we have
$$\limsup_{\mu\to\infty}\frac{1}{\mu}\log P(\mu A_x(\beta\mu,\mu)\in F)
\leq-\inf_{z\in F}xH_B(z/x;\beta,1,s_M)\ \mbox{for all closed sets}\ F$$
and
$$\liminf_{\mu\to\infty}\frac{1}{\mu}\log P\left(\mu A_{x}(\beta\mu,\mu)\in G\right)
\geq-\inf_{z\in G\cap\mathcal{E}}xH_B(z/x;\beta,1,s_M)\ \mbox{for all open sets}\ G$$
where $\mathcal{E}=(x\tilde{z}(\beta,1,s_M),\infty)$ is the set of exposed points of
$xH_B(\cdot/x;\beta,1,s_M)$.
\end{remark}

Now we study moderate deviations, i.e. with the analogue of Proposition 4.2 in \cite{MacciMartinucciPirozzi}.

\begin{proposition}\label{prop:MD2}
For every family of positive numbers $\{\varepsilon_\mu:\mu>0\}$ such that
\begin{equation}\label{eq:MD2-conditions}
\varepsilon_\mu\to 0\ \mbox{and}\ \mu\varepsilon_\mu\to\infty,
\end{equation}
the family $\left\{\sqrt{\mu\varepsilon_\mu}(A_x(\beta\mu,\mu)-\mathbb{E}[A_x(\beta\mu,\mu)]):\mu>0\right\}$
satisfies the LDP with speed $1/\varepsilon_\mu$, and good rate function $\tilde{I}_2$ defined by
$\tilde{I}_2(z):=\frac{z^2}{2x\Lambda^{\prime\prime}(0;\beta,1)}$, 
where $\Lambda^{\prime\prime}(0;\beta,1)=\frac{8\beta}{(\beta-1)^3}$.
\end{proposition}
\begin{proof}
We follow the lines of the proof of Proposition 4.2 in \cite{MacciMartinucciPirozzi} and we have to prove that
$$\lim_{x\to\infty}\frac{1}{1/\varepsilon_\mu}
\log\mathbb{E}\left[e^{\frac{s}{\varepsilon_\mu}\sqrt{\mu\varepsilon_\mu}(A_x(\beta\mu,\mu)-\mathbb{E}[A_x(\beta\mu,\mu)])}\right]
=\frac{x\Lambda^{\prime\prime}(0;\beta,1)}{2}s^2\ (\mbox{for all}\ s\in\mathbb{R}).$$
Firstly we note that
\begin{multline*}
\frac{1}{1/\varepsilon_\mu}
\log\mathbb{E}\left[e^{\frac{s}{\varepsilon_\mu}\sqrt{\mu\varepsilon_\mu}(A_x(\beta\mu,\mu)-\mathbb{E}[A_x(\beta\mu,\mu)])}\right]
=\varepsilon_\mu\left(\log G_{A_x(\beta\mu,\mu)}\left(s\sqrt{\frac{\mu}{\varepsilon_\mu}}\right)
-\mathbb{E}[A_x(\beta\mu,\mu)]s\sqrt{\frac{\mu}{\varepsilon_\mu}}\right)\\
=\mu\varepsilon_\mu\left(\frac{1}{\mu}\log G_{A_x(\beta\mu,\mu)}\left(\frac{\mu s}{\sqrt{\mu\varepsilon_\mu}}\right)-
\mathbb{E}[A_x(\beta\mu,\mu)]\frac{s}{\sqrt{\mu\varepsilon_\mu}}\right),
\end{multline*}
where $\frac{s}{\sqrt{\mu\varepsilon_\mu}}$ is close to zero if $\mu$ is large enough. Moreover
\begin{equation}\label{eq:expected-value-2}
\mathbb{E}[A_x(\beta\mu,\mu)]=[\log G_{A_x(\beta\mu,\mu)}]^\prime(0)
=\frac{2}{\mu(\beta-1)}+x\underbrace{\frac{\beta+1}{\beta-1}}_{=\Lambda^\prime(0;\beta,1)}
\end{equation}
(the second equality can be checked with some computations). Then, by Proposition \ref{prop:mgf}, 
eq. \eqref{eq:expected-value-2} and the Taylor formula of order 2 for the function $\Lambda(s;\beta,1)$, for $\mu$ large enough we have
\begin{multline*}
\frac{1}{1/\varepsilon_\mu}
\log\mathbb{E}\left[e^{\frac{s}{\varepsilon_\mu}\sqrt{\mu\varepsilon_\mu}(A_x(\beta\mu,\mu)-\mathbb{E}[A_x(\beta\mu,\mu)])}\right]\\
=\mu\varepsilon_\mu\left(\frac{1}{\mu}\log G_M\left(\log G_{C_0(\beta\mu,\mu)}\left(\frac{\mu s}{\sqrt{\mu\varepsilon_\mu}}\right)\right)
+\frac{x}{\mu}\Lambda\left(\frac{\mu s}{\sqrt{\mu\varepsilon_\mu}};\beta\mu,\mu\right)\right.\\
\left.-\left(\frac{2}{\mu(\beta-1)}+x\Lambda^\prime(0;\beta,1)\right)\frac{s}{\sqrt{\mu\varepsilon_\mu}}\right)\\
=\mu\varepsilon_\mu\left(\frac{1}{\mu}\log G_M\left(\log G_{C_0(\beta,1)}\left(\frac{s}{\sqrt{\mu\varepsilon_\mu}}\right)\right)
+x\Lambda\left(\frac{s}{\sqrt{\mu\varepsilon_\mu}};\beta,1\right)\right.\\
\left.-\left(\frac{2}{\mu(\beta-1)}+x\Lambda^\prime(0;\beta,1)\right)\frac{s}{\sqrt{\mu\varepsilon_\mu}}\right)\\
=\varepsilon_\mu\left(\log G_M\left(\log G_{C_0(\beta,1)}\left(\frac{s}{\sqrt{\mu\varepsilon_\mu}}\right)\right)
-\frac{2}{\beta-1}\frac{s}{\sqrt{\mu\varepsilon_\mu}}\right)
+x\mu\varepsilon_\mu\left(\frac{\Lambda^{\prime\prime}(0;\beta,1)}{2}\frac{s^2}{\mu\varepsilon_\mu}+o\left(\frac{1}{\mu\varepsilon_\mu}\right)\right)
\end{multline*}
and we conclude by taking the limit as $\mu\to\infty$.
\end{proof}

In the final part we present a non-central moderate deviation result. We start with the analogue of Lemma 4.1 in 
\cite{MacciMartinucciPirozzi}. 

\begin{lemma}\label{lem:weak-convergence}
For $\beta>1$, the random variables $\left\{\mu A_{x/\mu}(\beta\mu,\mu):\mu>0\right\}$ are equally distributed.
\end{lemma}
\begin{proof}
The result can be easily proved by taking the moment generating functions of the involved random variables, and by 
referring to the formulas presented in Proposition \ref{prop:mgf} (note that $G_{C_0(\beta\mu,\mu)}(\mu s)=G_{C_0(\beta,1)}(s)$
and $\Lambda(\mu s;\beta\mu,\mu)=\mu\Lambda(s;\beta,1)$). In fact these moment generating functions do not depend on $\mu$.
We omit the details.
\end{proof}

Now we prove the analogue of Proposition 4.3 in \cite{MacciMartinucciPirozzi}.
 
\begin{proposition}\label{prop:MD2-nc}
Assume that $s_M\geq\log\sqrt{\beta}$, for $\beta>1$. Then, for every family of positive numbers $\{\varepsilon_\mu:\mu>0\}$ 
such that \eqref{eq:MD2-conditions} holds, the family $\left\{\mu\varepsilon_\mu A_{x/(\mu\varepsilon_\mu)}(\beta\mu,\mu):\mu>0\right\}$
satisfies the LDP with speed $1/\varepsilon_\mu$, and good rate function $I_2$ (see Proposition \ref{prop:LD2}).
\end{proposition}
\begin{proof}
We follow the lines of the proof of Proposition 4.3 in \cite{MacciMartinucciPirozzi}. We consider Proposition \ref{prop:mgf}
(again we note that $G_{C_0(\beta\mu,\mu)}(\mu s)=G_{C_0(\beta,1)}(s)$ and $\Lambda(\mu s;\beta\mu,\mu)=\mu\Lambda(s;\beta,1)$),
Remark \ref{rem:the-domain-of-mgf-of-absorption-time} (Case A) and Lemma \ref{lem:Legendre}(i). We distiguish two cases.

In the first case, i.e. if $s_M=\log\sqrt{\beta}$ and $\mathcal{D}(G_M)$ is open, then
\begin{multline*}
\frac{1}{1/\varepsilon_\mu}\log\mathbb{E}\left[e^{\frac{s}{\varepsilon_\mu}\mu\varepsilon_\mu A_{x/(\mu\varepsilon_\mu)}(\beta\mu,\mu)}\right]
=\varepsilon_\mu\log\mathbb{E}\left[e^{s\mu A_{x/(\mu\varepsilon_\mu)}(\beta\mu,\mu)}\right]\\
=\left\{\begin{array}{ll}
\varepsilon_\mu\{\log G_M(\log G_{C_0(\beta\mu,\mu)}(s\mu))+\frac{x}{\mu\varepsilon_\mu}\Lambda(s\mu;\beta\mu,\mu)\}
&\ \mbox{for}\ s\mu<\frac{(\sqrt{\beta\mu}-\sqrt{\mu})^2}{2}\\
\infty&\ \mbox{otherwise}
\end{array}\right.\\
=\left\{\begin{array}{ll}
\varepsilon_\mu\log G_M(\log G_{C_0(\beta,1)}(s))+x\Lambda(s;\beta,1)
&\ \mbox{for}\ s<\frac{(\sqrt{\beta}-1)^2}{2}\\
\infty&\ \mbox{otherwise},
\end{array}\right.
\end{multline*}
and therefore
$$\lim_{\mu\to\infty}\frac{1}{1/\varepsilon_\mu}\log\mathbb{E}\left[e^{\frac{s}{\varepsilon_\mu}\mu\varepsilon_\mu
A_{x/(\mu\varepsilon_\mu)}(\beta\mu,\mu)}\right]=\left\{\begin{array}{ll}
x\Lambda(s;\beta,1)&\ \mbox{if}\ s<\frac{(\sqrt{\beta}-1)^2}{2}\\
\infty&\ \mbox{if}\ s\geq\frac{(\sqrt{\beta}-1)^2}{2}.
\end{array}\right.$$

In the second case, i.e. if $s_M>\log\sqrt{\beta}$ and/or $\mathcal{D}(G_M)$ is closed, then
$$\frac{1}{1/\varepsilon_\mu}\log\mathbb{E}\left[e^{\frac{s}{\varepsilon_\mu}\mu\varepsilon_\mu A_{x/(\mu\varepsilon_\mu)}(\beta\mu,\mu)}\right]
=\left\{\begin{array}{ll}
\varepsilon_\mu\log G_M(\log G_{C_0(\beta,1)}(s))+x\Lambda(s;\beta,1)
&\ \mbox{for}\ s\leq\frac{(\sqrt{\beta}-1)^2}{2}\\
\infty&\ \mbox{otherwise},
\end{array}\right.$$
and therefore
$$\lim_{\mu\to\infty}\frac{1}{1/\varepsilon_\mu}\log\mathbb{E}\left[e^{\frac{s}{\varepsilon_\mu}\mu\varepsilon_\mu
A_{x/(\mu\varepsilon_\mu)}(\beta\mu,\mu)}\right]=\left\{\begin{array}{ll}
x\Lambda(s;\beta,1)&\ \mbox{if}\ s\leq\frac{(\sqrt{\beta}-1)^2}{2}\\
\infty&\ \mbox{if}\ s>\frac{(\sqrt{\beta}-1)^2}{2}.
\end{array}\right.$$

Finally, in both cases, the desired LDP (for every choice of positive numbers $\{\varepsilon_\mu:\mu>0\}$ such that 
\eqref{eq:MD2-conditions} holds) can be obtained as a straightforward application of Theorem \ref{th:GE}.
\end{proof}

\section{Numerical estimates by simulations}\label{sec:numerical-results}
In this section we follow the same lines of Section 5 in \cite{MacciMartinucciPirozzi}. We refer to an asymptotic Normality 
result under the scaling 2, i.e. the weak convergence of $\sqrt{\mu}(A_x(\beta\mu,\mu)-\mathbb{E}[A_x(\beta\mu,\mu)])$ to 
the centered Normal distribution with variance $x\Lambda^{\prime\prime}(0;\beta,1)$. The aim is to present some numerical 
values obtained by simulations to estimate $\beta$; actually we assume that $\beta>\beta_0$ for some known $\beta_0>1$.

Now we recall some formulas presented in Section 5 in \cite{MacciMartinucciPirozzi}. Let $\Phi$ be the standard Normal 
distribution function. We denote the simulated sample mean of $\overline{A}_x(\beta\mu,\mu)$ for chosen values 
$\beta=\beta_*>\beta_0>1$ by $\overline{A}_x(\beta_*\mu,\mu)$ and, when $\mu$ is large, we have:
\begin{itemize}
\item the confidence interval for $\beta_*$ at the level $\ell$, when $x<\overline{A}_x(\beta_*\mu,\mu)-
\sqrt{\frac{8\beta_0x}{(\beta_0-1)^3}}\frac{\Phi^{-1}\left(\frac{1+\ell}{2}\right)}{\sqrt{\mu}}$,
\begin{equation}\label{eq:confidence-interval-for-beta}
\left(\frac{\overline{A}_x(\beta_*\mu,\mu)+\sqrt{\frac{8\beta_0x}{(\beta_0-1)^3}}\frac{\Phi^{-1}\left(\frac{1+\ell}{2}\right)}{\sqrt{\mu}}+x}
{\overline{A}_x(\beta_*\mu,\mu)+\sqrt{\frac{8\beta_0x}{(\beta_0-1)^3}}\frac{\Phi^{-1}\left(\frac{1+\ell}{2}\right)}{\sqrt{\mu}}-x},
\frac{\overline{A}_x(\beta_*\mu,\mu)-\sqrt{\frac{8\beta_0x}{(\beta_0-1)^3}}\frac{\Phi^{-1}\left(\frac{1+\ell}{2}\right)}{\sqrt{\mu}}+x}
{\overline{A}_x(\beta_*\mu,\mu)-\sqrt{\frac{8\beta_0x}{(\beta_0-1)^3}}\frac{\Phi^{-1}\left(\frac{1+\ell}{2}\right)}{\sqrt{\mu}}-x}\right);
\end{equation}
\item the point estimation of $\beta_*$
\begin{equation}\label{eq:point-estimation-for-beta}
\frac{\overline{A}_x(\beta_*\mu,\mu)+x}{\overline{A}_x(\beta_*\mu,\mu)-x}.
\end{equation}
\end{itemize}

Now we are ready to present some numerical values for Example \ref{ex:shifted-Poisson-distribution} (instead of Example 
\ref{ex:geometric-distribution} as in \cite{MacciMartinucciPirozzi}). In all cases we perform simulations by setting $x=1$ and 
$\beta_0=1.50$; furthermore, the size of simulated sample paths is $10^3$ and the confidence level is $\ell=0.95$. For
each table we vary one parameter (among $\theta$, $\mu$ and $\beta_*$) and the other two are fixed.

\begin{table}[h!]
	\begin{center}
		\caption{Numerical approximations for the confidence interval for $\beta$ varying $\theta$}
		\label{tab:table1}
		\begin{tabular}{||c||c|c|c|c|c|c} 
			\hline
			$\theta$ & $\mu$ &$\beta_*$ & $x\frac{\beta_*+1}{\beta_*-1}$&$\overline{A}_x(\beta_* \mu, \mu)$
			&Confidence Interval \eqref{eq:confidence-interval-for-beta}& Point Estimation \eqref{eq:point-estimation-for-beta}\\
			\hline
			1.5 & 1000 & 1.75 & 3.666667 & 3.672570 &(1.657130,1.868960) & 1.748343\\
			\hline
			3 & 1000 & 1.75 & 3.666667 & 3.661619 & (1.65950,1.873114) & 1.751422\\
			\hline
			5 & 1000 & 1.75 & 3.666667 & 3.672501 & (1.657145,1.868985) & 1.748363\\
			\hline
			10 & 1000 & 1.75 & 3.666667 & 3.698365 & (1.651608,1.859329) & 1.741190\\
			\hline
		\end{tabular}
	\end{center}
\end{table}

\begin{table}[h!]
	\begin{center}
		\caption{Numerical approximations for the confidence interval for $\beta$ varying $\mu$}
		\label{tab:table2}
		\begin{tabular}{c||c||c|c|c|c|c} 
			\hline
			$\theta$ & $\mu$ &$\beta_*$ & $x\frac{\beta_*+1}{\beta_*-1}$&$\overline{A}_x(\beta_* \mu, \mu)$
			&Confidence Interval \eqref{eq:confidence-interval-for-beta}& Point Estimation \eqref{eq:point-estimation-for-beta}\\
			\hline
			3 & 1000 & 2 & 3 & 3.00749 & (1.840883,2.222106) & 1.996271\\
			\hline
			3 & 5000 & 2 & 3 & 2.999793 & (1.9234911,2.09058) & 2.000104\\
			\hline
			3 & 10000 & 2 & 3 & 2.999903 & (1.944638,2.062365) & 2.000048\\
			\hline
			3  & 50000 & 2 & 3 & 2.999883 & (1.974494,2.026999) & 2.000058\\
			\hline
		\end{tabular}
	\end{center}
\end{table}

\begin{table}[h!]
	\begin{center}
		\caption{Numerical approximations for the confidence interval for $\beta$ varying $\beta_*$}
		\label{tab:table3}
		\begin{tabular}{c|c||c||c|c|c|c} 
			\hline
			$\theta$ & $\mu$ &$\beta_*$ & $x\frac{\beta_*+1}{\beta_*-1}$&$\overline{A}_x(\beta_* \mu, \mu)$
			&Confidence Interval \eqref{eq:confidence-interval-for-beta}& Point Estimation \eqref{eq:point-estimation-for-beta}\\
			\hline
			5 & 1000 & 1.5 & 5 & 5.018858 & (1.455599,1.548262) & 1.497654\\
			\hline
			5 & 1000 & 2 & 3 & 3.010646 & (1.839767,2.21971) & 1.994705\\
			\hline
			5 & 1000 & 2.5 & 2.${\overline 3}$ & 2.338547 & (2.169924,3.067012) & 2.494158\\
			\hline
			5 & 1000 & 3 & 2 & 2.000227 & (2.458583,4.178333) & 2.999545\\
			\hline
		\end{tabular}
	\end{center}
\end{table}

In Table \ref{tab:table1} we consider some values of $\theta$ ($\mu$ and $\beta_*$ are constant). The point estimates 
and the length of the confidence intervals are stable for $\theta=1.5$, $\theta=3$ and $\theta=5$; on the contrary, for 
$\theta=10$, the point estimate is slightly less accurate, and the confidence interval is slightly narrower.

In Table \ref{tab:table2} we consider some large values of $\mu$ ($\theta$ and $\beta_*$ are constant). Then, as one
can expect, the accuracy of point estimates and confidence intervals generally improves when $\mu$ increases; indeed
our formulas \eqref{eq:confidence-interval-for-beta} and \eqref{eq:point-estimation-for-beta} concern the scaling 2 
where $\mu\to\infty$.

In Table \ref{tab:table3} we consider some values of $\beta_*$ ($\theta$ and $\mu$ are constant). In this case, as 
$\beta_*$ increases, we have more accurate point estimates and wider confidence intervals.

\paragraph{Funding.} Antonella Iuliano acknowledges the partial support of Indam-GNCS (project ``Mo\-delli di shock basati sul processo di 
conteggio geometrico e applicazioni alla sopravvivenza'' (CUP E55F22000270001)).
Claudio Macci acknowledges the partial support of MIUR Excellence Department Project awarded to the Department of Mathematics, University 
of Rome Tor Vergata (CUP E83C18000100006), of University of Rome Tor Vergata (project ``Asymptotic Methods in Probability'' 
(CUP E89C20000680005) and project ``Asymptotic Properties in Probability'' (CUP E83C22001780005)) and of Indam-GNAMPA.


\end{document}